\documentclass[twocolumn,prb,showpacs]{revtex4-1}
\usepackage{graphicx}
\usepackage{epstopdf}
\usepackage{dcolumn}
\usepackage{amsmath}
\usepackage{amsfonts}
\usepackage{multirow} 
\usepackage[active]{srcltx}
\usepackage{hyperref}
\usepackage{yfonts}
\usepackage{hyperref}

\begin{document}

\title{Solution of the Mayan Calendar Enigma}

\author{Thomas Chanier$^{1*}$}

\affiliation{$^{1}$Independent researcher, 1025 12$^\mathrm{th}$ avenue, Coralville, Iowa 52241, USA}

\begin{abstract}
The Mayan calendar is proposed to derive from an arithmetical model of naked-eye astronomy. The Palenque and Copan lunar equations, used during the Maya Classic period (200 to 900 AD) are solution of the model and the results are expressed as a function of the Xultun numbers, four enigmatic Long Count numbers deciphered in the Maya ruins of Xultun, dating from the IX century AD, providing strong arguments in favor of the use of the model by the Maya. The different Mayan Calendar cycles can be derived from this model and the position of the Calendar Round at the mythical date of creation 13(0).0.0.0.0 4 Ahau 8 Cumku is calculated. This study shows the high proficiency of Mayan mathematics as applied to astronomy and timekeeping for divinatory purposes.\footnote{Published in the Mathematical Intelligencer: A Possible Solution to the Mayan Calendar Enigma, Math. Intell. $\mathbf{40}$, 18-25 (2018).} \end{abstract}


\maketitle








\section{Introduction}

Mayan priests-astronomers were known for their astronomical and mathematical proficiency, as demonstrated in the Dresden Codex, a XIV century AD bark-paper book containing accurate astronomical almanacs aiming to correlate ritual practices with astronomical observations. However, due to the zealous role of the Inquisition during the XVI century AD Spanish conquest of Mexico, most of these codices were destroyed and only four of them, the Dresden Codex, the Madrid Codex, the Paris Codex and the Grolier Codex remain today the only original native written records on Maya civilization. Thanks to the work of Mayan archeologists and epigraphists since the early XX century, these four codices, along with numerous inscriptions on monuments, were deciphered, underlying the Mayan cyclical concept of time. This is demonstrated by the Mayan Calendar formed by a set of three interlocking cycles: the Calendar Round, the Kawil-directions-colors cycle and the Long Count Calendar.

\begin{figure}[h!]
\begin{center}
\includegraphics[width = 75mm,angle=0]{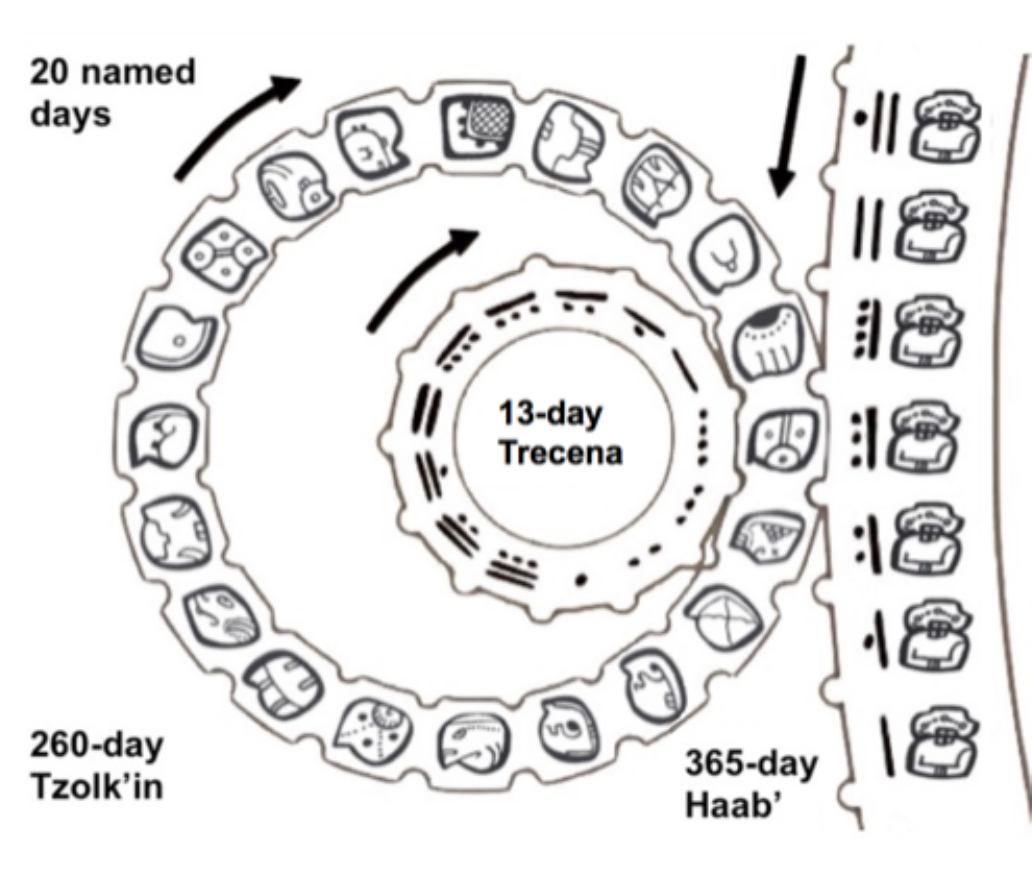}
\caption{Contemporary representation of the Calendar Round. The Tzolk'in day and Haab' month glyphs are taken from Ref. [\onlinecite{Bricker2011_68}] and the Mayan vigesimal system is used (a dot represents 1 and a bar 5). The 260-day Tzolk'in, obtained by the permutation of 20 named days and the 13-day Trecena, coincides with the 365-day Haab' every 18980 days = 52 $\times$ 365 = 73 $\times$ 260 corresponding to the Calendar Round. The setting corresponds to the day 4 Ahau 8 Cumku, day origin of the Calendar Round at the mythical date of creation. The next day is 5 Imix 9 Cumku.}\label{CR}
\end{center}
\end{figure}

The Calendar Round (CR) represents a day in a 18980-day cycle, a period of roughly 52 years, the combination of the 365-day solar year Haab' and the 260-day religious year Tzolk'in. The Tzolk'in combines the 13-day Trecena cycle (numerated from 1 to 13) with 20 named days (Imix, Ik, Akbal, Kan, Chicchan, Cimi, Manik, Lamat, Muluc, Oc, Chuen, Eb, Ben, Ix, Men, Cib, Caban, Etznab, Cauac, and Ahau). This forms a list of 260 ordered Tzolk'in dates from 1 Imix, 2 Ik, ... to 13 Ahau [\onlinecite{Aveni2001_143}]. The Haab' comprises 18 named months (Pop, Uo, Zip, Zotz, Tzec, Xul, Yaxkin, Mol, Chen, Yax, Zac, Ceh, Mac, Kankin, Muan, Pax, Kayab, and Cumku) with 20 days each (Winal) plus 1 extra month (Uayeb) with 5 nameless days. This forms a list of 365 ordered Haab' dates from 0 Pop, 1 Pop, ... to 4 Uayeb [\onlinecite{Aveni2001_147}]. The Tzolk'in and the Haab' coincide every 73 Tzolk'in or 52 Haab' or a Calendar Round, the least common multiple (LCM) 1 CR = LCM(260,365) = $73\times260=52\times365=18980$ days. In the Calendar Round, a date is represented by $\alpha$X$\beta$Y with the religious month $1\leq\alpha\leq13$, X one of the 20 religious days, the civil day $0\leq\beta\leq19$, and Y one of the 18 civil months, $0\leq\beta\leq4$ for the Uayeb. Fig. \ref{CR} shows a contemporary representation of the Calendar Round as a set of three interlocking wheels: the Tzolk'in, formed by a 13-month and a 20-day wheels and the Haab'.

\begin{figure}[h!]
\begin{center}
\includegraphics[width = 55mm,angle=0]{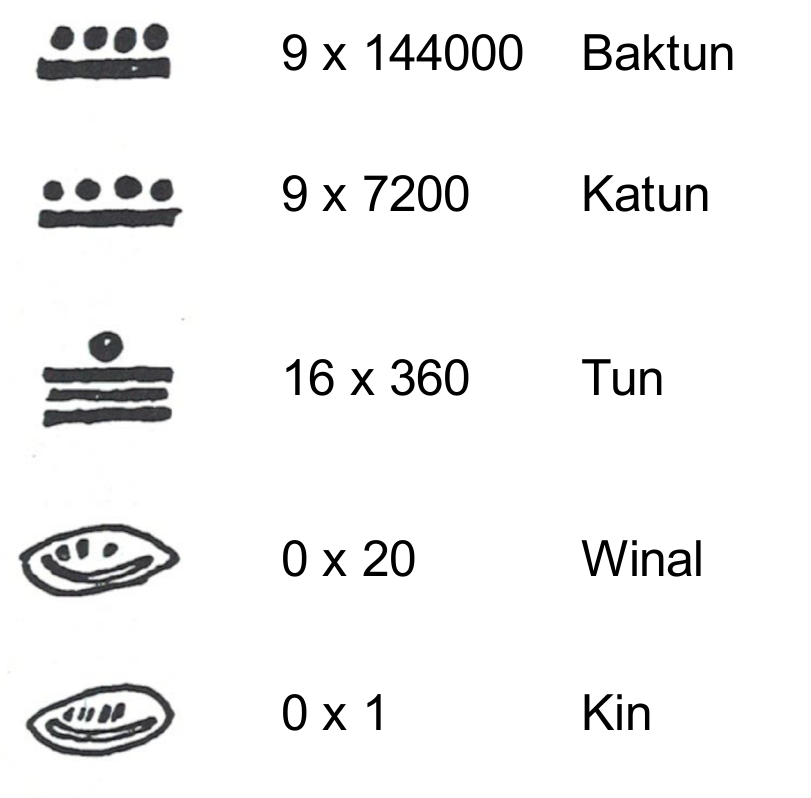}\\[.3cm]
\caption{Long Round number $\cal{LR}$ = 9.9.16.0.0 = 9 $\times$ 144000 + 9 $\times$ 7200 + 16 $\times$ 360 + 0 $\times$ 20 + 0 $\times$ 1 = 1366560 present on page 24 of the Dresden Codex. Long Count numbers are represented vertically using the Mayan vigesimal system: a shell represents 0, a dot 1 and a bar 5.}\label{LR}
\end{center}
\end{figure}

For longer periods of time, the Maya used the Long Count Calendar (LCC), describing a date $D$ in a 1872000-day Maya Era of 13 Baktun, a religious cycle of roughly 5125 years, counting the number of day elapsed since the Mayan origin of time. This mythical date of creation, the day 13(0).0.0.0.0 4 Ahau 8 Cumku, is carved on Stela C (775 AD) of the Maya site Quirigua (present-day Guatemala)  [\onlinecite{VanStone,Quirigua}] and corresponds to the Gregorian Calendar date 11 August 3114 BC according to the Goodman-Martinez-Thompson (GMT) correlation [\onlinecite{Aveni2001_136,Bricker2011_71}]. An interesting example of Long Count number can be found in the introduction of the Venus table on page 24 of the Dresden Codex: the so-called Long Round number $\cal{LR}$ = 9.9.16.0.0 = 1366560 days, a whole multiple of the Tzolk'in, the Haab', the Calendar Round, the Tun, Venus and Mars synodic periods: $\cal{LR}$ = 5256 $\times$ 260 = 3744 $\times$ 365 = 72 $\times$ 18980 = 3796 $\times$ 360 = 2340 $\times$ 584 = 1752 $\times$ 780. The Long Round number (Fig. \ref{LR}) can be expressed as $\cal{LR}$ = 9 $\times$ 144000 + 9 $\times$ 7200 + 16 $\times$ 360 + 0 $\times$ 20 + 0 $\times$ 1 as a function of the Long Count periods (the 1-day Kin, the 20-day Winal, the 360-day Tun, the 7200-day Katun and the 144000-day Baktun) [\onlinecite{Aveni2001_191}]. The Long Count periods are commensurate with the Tzolk'in and the Haab': \{LCM(260,$P_i$)/$P_i$ = 13, LCM(365,$P_i$)/$P_i$ = 73, $P_i=18\times20^i,\ i>0$\}. The XXI century saw the passage of a new Maya Era on the day of the winter solstice 21 December 2012 (GMT correlation) or 13(0).0.0.0.0 4 Ahau 3 Kankin, a date carved on Monument 6 of Tortuguero (present-day Mexico), a Maya stone from the VII century AD [\onlinecite{Stuart2012_25}].

The Kawil-directions-colors cycle or 4-Kawil is a 3276-day cycle, the combination of the 819-day Kawil and the 4 directions-colors (East-Red, South-Yellow, West-Black, North-White) [\onlinecite{Berlin1961}]. The 4-Kawil counts the number of day in four 819-day months (each of them corresponding to one direction-color) in a non-repeating 3276-day cycle. At the mythical date of creation, the Kawil count is 3 and the direction-color is East-Red. A Kawil date is then defined as $D\equiv$  mod($D$ + 3, 819) and the direction-color is given by $n$ = mod(int(($D$+3)/819), 4), $n$ = \{0, 1, 2 ,3\} = \{East-Red, South-Yellow, West-Black, North-White\}. Although several myths exist around Mayan religion, the origin of the different Mayan Calendar cycles remains unknown.

A complete Maya date contains also a glyph $G_i$ with $i$ = $1\dots9$ corresponding to the 9 Lords of the Night and the lunar series: the 29(30)-day Moon age (number of days elapsed in the current lunation) and a lunation count (number of lunation in a series of five or six). The calculation of the Moon age in the lunar series of a mythical date $LC$ is calculated from the new Moon date $LC_0$ as: $MA$ = remainder of ($LC$ - $LC_0$)/$S$ where $S=n/m$ is the Moon ratio corresponding to the lunar equation $m$ lunations = $n$ days [\onlinecite{Fuls2007}]. Mayan priests-astronomers used particular lunar equations such as 149 lunations = 4400 days (Copan Moon ratio) and 81 lunations = 2392 days (Palenque formula). The Palenque formula corresponds to a Moon synodic period of 29.530864 days, differing by only 24 seconds from the modern value (29.530588 days) [\onlinecite{Fuls2007,Aveni2001_163,Teeple1930}]. It is unclear how the Maya determined these particular lunar equations.

\begin{table}[h!]
\begin{center}
\begin{tabular}{l l l c l l}
$\mathcal{X}_i$ & LCC & $D$ [day] & $\mathcal{X}_i/56940$ & \\\hline
$\mathcal{X}_0$ &	2.7.9.0.0 & 341640 & 6\\
$\mathcal{X}_1$ &	8.6.1.9.0 & 1195740 & 21\\
$\mathcal{X}_2$ &	12.5.3.3.0 & 1765140 & 31\\
$\mathcal{X}_3$ &	17.0.1.3.0 & 2448420 & 43\\
\hline
\end{tabular}
\caption{Xultun numbers $\mathcal{X}_i$ [\onlinecite{Saturno2012}]. 56940 = LCM(365,780) is their largest common divisor and their sum is such that $\sum_{0}^3\mathcal{X}_i$ = 101 $\times$ 56940.}\label{Xultun}
\end{center}
\end{table}

In 2012, four Long Count numbers, the Xultun numbers (Table \ref{Xultun}) and three lunar tables, have been discovered on the walls of a small painted room in the Maya ruins of Xultun (present-day Guatemala), dating from the early IX century AD [\onlinecite{Saturno2012}]. These numbers have a potential astronomical meaning. Indeed, $\mathcal{X}_0$ is a whole multiple of Venus and Mars synodic periods: 341640 =  585 $\times$ 584 = 438 $\times$ 780. $\mathcal{X}_0$ = $\cal{LR}$/4 = LCM(260,360,365) is the commensuration of the Tzolk'in, the Tun and the Haab' and $\mathcal{X}_1=365\times3276$ is the commensuration of the Haab' and the 4-Kawil. However, the meaning of $\mathcal{X}_2$ and $\mathcal{X}_3$ is unknown. The greatest common divisor of the $\mathcal{X}_i$'s is 56940 = LCM(365,780) = 3 CR, the commensuration of the Haab' and Mars synodic period. The three Xultun lunar tables, corresponding to a time span of 4429 (12.5.9), 4606 (12.14.6) and 4784 (13.5.4) days were attributed to solar/lunar eclipse cycles due to similarities in structure with the Dresden Codex eclipse table [\onlinecite{Saturno2012}]. It was noted that 4784 = 2 $\times$ 2392 days represented 162 lunations, twice that of the Palenque lunar reckoning system 81 lunations = 2392 days [\onlinecite{Teeple1930}]. The length of the Dresden Codex eclipse table 11960 = 5 $\times$ 2392 days = 405 lunations corresponds to five times the Palenque formula [\onlinecite{Bricker1983}]. The lengths of the solar/lunar eclipse tables are unexplained.\\[.3cm]

The paper is organized as follows. In section 2, a model of naked-eye astronomy is described. The Palenque and Copan lunar equations are calculated from the calendar super-number defined as the least common multiple of 9 astronomical input parameters describing the synodic movements of the Moon and the five planets visible in the night sky with a naked-eye (Mercury, Venus, Mars, Jupiter and Saturn). This gives a first clue that the model can be attributed to the Maya. In section 3, the Mayan Calendar, combination of the Calendar Round, the Long Count Calendar and the Kawil-directions-colors cycle is derived from the model. The results are expressed as a function of the Xultun numbers, providing an additional argument in favor of the model. The initialization of the Calendar Round at the Mayan mythical date of creation is calculated. In section 4, important Mayan mythical dates are deduced from the model and their numerological significance is discussed in terms of culturally and historically important events. Section 5 is left for conclusion.

\section{Mayan model of naked-eye astronomy}

The level of sophistication displayed in the Dresden Codex suggests the high astronomical proficiency of Mayan priests-astronomers. It is therefore reasonable to assume that the Maya measured the synodic periods of the five planets Mercury, Venus, Mars, Jupiter and Saturn visible by naked-eye observation of the night sky. Their canonic synodic periods are given in Table \ref{SP}. Evidence of their use has been found in different Mayan codices for Mercury, Venus and Mars [\onlinecite{Bricker2011_163}]. References to Jupiter and Saturn have been found in some Maya texts [\onlinecite{Fox,Lounsbury1989,Aveni1994}] and their retrograde motions were associated to the Maya Katun cycle [\onlinecite{Milbrath}]. The three relevant lunar months are the two lunar semesters of 177 and 178 days (6 Moon synodic periods) and the pentalunex of 148 days (5 Moon synodic periods), parameters used for the prediction of solar/lunar eclipses in the Dresden Codex eclipse table [\onlinecite{Bricker1983}]. From the prime factorizations of the 9 astronomical input parameters (Table \ref{SP}), we calculate the calendar super-number $\mathcal{N}$ defined as the least common multiple of the $P_i$'s:
\begin{eqnarray}\label{astro}
\mathcal{N} & = & 768039133778280\\\nonumber
 & = & 2^3 \times 3^3\times 5\times 7 \times 13 \times 19 \times 29 \times 37\\\nonumber
&\ &\times\ 59 \times 73 \times 89\\\nonumber
 & = & 365 \times 3276 \times 2 \times 3 \times 19 \times 29 \times 37\\\nonumber
&\ &\times\ 59 \times 89\\\nonumber
 & = & \mathrm{LCM}(360,365,3276) \times 3 \times 19 \times 29\\\nonumber &\ & \times\ 37 \times 59 \times 89
\end{eqnarray}

\begin{table}[b!]
\begin{center}
\begin{tabular}{l l c l}
Planet &$i$ &  $P_i$ [day] & Prime factorization\\\hline
Mercury &1 & 116 & 2$^2$ $\times$ 29\\
Venus &2 & 584 & 2$^3$ $\times$ 73\\
Earth &3 & 365 & 5 $\times$ 73\\
Mars &4 & 780 & 2$^2$ $\times$ 3 $\times$ 5 $\times$ 13\\
Jupiter &5 & 399 & 3 $\times$ 7 $\times$ 19\\
Saturn &6 & 378 & 2 $\times$ 3$^3$ $\times$ 7\\\hline
Lunar &7 & 177 & 3 $\times$ 59\\
senesters &8 &  178 & 2 $\times$ 89\\
Pentalunex & 9 &  148 & 2$^2$ $\times$ 37\\\hline
\end{tabular}
\caption{Prime factorization of the planet canonic synodic periods and the three Mayan lunar months [\onlinecite{Bricker1983}].}\label{SP}
\end{center}
\end{table}

Equ. \ref{astro} gives the calendar super-number and its prime factorization. It is expressed as a function of the Tun (360 = 18 $\times$ 20), the Haab' (365 = 5 $\times$ 73) and the 4-Kawil (3276 = 2$^2$ $\times$ $3^2$ $\times$ 7 $\times$ 13). The solar year Haab' and the $P_i$'s are relatively primes (except Venus and Mars): the \{LCM($P_i$,365)/365, $i$ = $1\dots9$\} = \{116, 8, 1, 156, 399, 378, 177, 178, 148\} (Table \ref{SP}). The 4-Kawil is defined as the \{LCM($P_i$,3276)/3276, $i$ = $1\dots9$\} = \{29, 146, 365, 5, 19, 3, 59, 89, 37\}. The Haab' and the 4-Kawil are relatively primes: the LCM(365,3276) = 365 $\times$ 3276 = $\mathcal{X}_1$ such as \{LCM($P_i$,$\mathcal{X}_1$)/$\mathcal{X}_1$, $i$ = $1\dots9$\} = \{29, 2, 1, 1, 19, 3, 59, 89, 37\}. 360 is the nearest integer to 365 such that the LCM(360,3276) = 32760 and the \{LCM($P_i$,32760)/32760, $i$ = $1\dots9$\} = \{29, 73, 73, 1, 19, 3, 59, 89\}. The Tun-Haab'-Kawil cycle is given by $\mathcal{Y}$ = LCM(360,365,3276) = 7 $\times$ $\mathcal{X}_0$ = 2391480 such as \{LCM($P_i$,$\mathcal{Y}$)/$\mathcal{Y}$, $i$ = $1\dots9$\} = \{29, 1, 1, 1, 19, 3, 59, 89, 37\}. 

\begin{table} [h!]
\begin{center}
\begin{tabular}{l c l c}
$T$ [day] & $L$ &  $S$ [day] & $\varepsilon$ [day]\\\hline
11960$^a$ & 405 &  29.530864 & 1\\
4784$^b$ &162 &   29.530864 & 1\\
4606$^b$ & 156 &  29.525641 & 8\\
4429$^b$ & 150 &  29.526667 & 11\\
4400$^c$ & 149 & 29.530201 & 2\\
2392$^d$ &  81 &  29.530864 & 1\\\hline
\multicolumn{2}{l}{Modern value} & 29.530588& 4\\  
\end{tabular}
\caption{Mayan lunar period $S$ = $T$/$L$ calculated from the length $T$ of the lunar tables and the corresponding number of lunations $L$ = Rd($T$,29.53) as compared to the modern value of the Moon synodic period. The length of the lunar tables are taken from the Dresden Codex eclipse table$^a$ [\onlinecite{Bricker1983}], the Xultun lunar table$^b$ [\onlinecite{Saturno2012}], the Copan Moon ratio$^c$ [\onlinecite{Aveni2001_163}] and the Palenque formula.$^d$ [\onlinecite{Teeple1930}] The error $\varepsilon$ of the Moon ratio is calculated from Equ. \ref{lunarequ}.}\label{lunar}
\end{center}
\end{table} 

The calculation of the Moon age in the lunar series of mythical dates necessitated a precise value of the Moon ratio corresponding to a particular lunar equation. Mayan priests-astronomers carefully recorded the lunar equation $L$ lunations = $T$ days for extended periods of time [\onlinecite{Justeson}] but it is unclear how they determined the Palenque and Copan lunar equations. A possible method can be derived from the calendar super-number by astronomical observation of the Moon, giving an idea of how they may have proceeded. First, they correlated the synodic movement of the Moon (using the pentalunex and the two lunar semesters) with the solar year and the five planet synodic periods corresponding to the calendar super-number $\mathcal{N}$. They were looking for a Moon ratio $S=T/L$ such as $\mathcal{N}/S$ is an integer with the error $\varepsilon\rightarrow0$:  
\begin{equation}\label{lunarequ}
\varepsilon=|\mathcal{N}-\mathrm{Rd}(\mathcal{N}/S) \times S|
\end{equation}
where Rd() is the nearest integer round function. The results are given in Table \ref{lunar}. In Palenque (present-day Mexico), somewhere between the III century BC and the VIII century AD, a careful analysis of the lunar data allowed Mayan priests-astronomers to determine the Palenque formula 81 lunations = 2392 days (1-day error). In Copan (present-day Honduras), somewhere between the V and the IX centuries AD, similar attempts allowed them to determine the Copan Moon ratio 149 lunations = 4400 days (2-day error). In Xultun (present-day Guatemala), in the IX century AD, Mayan priests-astronomers tried to find other solutions by considering three different periods close to the Copan value: 150 lunations = 4429 days (11-day error), 156 lunations = 4606 (8-day error) and 162 lunations = 4784 days (1-day error). It seems that from this date, a unified lunar ratio, the Palenque formula, was used up to the XIV century AD as shown in the Dresden Codex eclipse table with a lunar equation 405 lunations = 11960 days such as $S_0$ = 11960/405 = 4784/162 = 2392/81 = 2$^3$ $\times$ 13 $\times$ 23/3$^4$ = 29.530864 days. Indeed, the Palenque formula ($\varepsilon$ = 1.28 day) constitutes a slight improvement compared to the Copan Moon ratio ($\varepsilon$ = 1.88 day). The Palenque formula corresponds to the equation 81 $\times$ $\mathcal{N}$ + 104 = 26008014145502 $\times$ 2392 and the error $\varepsilon$ = 104/81 = 1.28 (Equ. \ref{lunarequ}). To perform such tedious calculations, Mayan priests-astronomers may have used a counting device, such as an abacus [\onlinecite{Thompson1941,Thompson1950}]. Archaeological evidence from the X century AD attests the use of an abacus by the Aztec [\onlinecite{Sanchez1961}]. The Aztec counting device, the Nepohualtzitzin, consisted of a wooden frame with threaded strings of maize kernel and may have originated from an early Mayan abacus. They also benefited from the use of the Mayan vigesimal system. The efficiency of Mayan numerals for arithmetical calculations has been noted previously [\onlinecite{Bietenholz2013,FrenchAnderson1971}]. We can legitimately ask whether the Maya were able to handle such huge super-numbers ($\mathcal{N}$ $\sim$ 7.68 $\times$ 10$^{14}$). The use of tremendous Long Count numbers has been identified on various monuments. For example, the inscriptions on  Coba Stela 1 (present-day Mexico), dating from the VII century AD, represents a mythical date [13. $\dots$ 13.] 13.0.0.0.0 4 Ahau 8 Cumku (20 coefficients 13 including the Baktun 144000 = 18 $\times$ 20$^3$) corresponding to a decimal number of the order of 1 $\times$ 10$^{31}$ (taking into account higher order Long Count periods 18 $\times$ 20$^i$, $i>3$) [\onlinecite{Fuls2007}].

A question arises about the choice of the particular values of the Palenque and Copan lunar equations. To answer this question, we calculate the lunar equations $i$ lunations = $T_i^0$ days where $T_i^0$ = Rd($i$ $\times$ $S_M$) (nearest integer approximation), its commensuration with the Tzolk'in LCM(260,$T_i^0$), the Moon ratio $S_i$ = $T_i^0$/$i$ and the corresponding error $\varepsilon_i$ (Equ. \ref{lunarequ}) for $i$ = 1 to 643 lunations ($T_{643}^0$ = 18988 $>$ 1 CR), taking into account the modern value of the Moon synodic period ($S_M$ = 29.530588 days). We consider the values minimizing  Equ. \ref{lunarequ} such as the Moon ratio $S_i$ = $T_i^0$/$i$ = LCM(260,$T_i^0$)/$j$ where $j$ = Rd(LCM(260,$T_i^0$)/$S_M$) is the number of lunations corresponding to the Tzolk'in-lunar commensuration period. The list contains the Palenque formula 2392/81 = 4784/162 = LCM(260,2392)/405 = 11960/405 = 29.530864 ($\varepsilon=1.28$) and the Copan Moon ratio 4400/149 = LCM(260,4400)/1937 = 57200/1937 = 29.530201 ($\varepsilon$ = 1.88 day). The Palenque formula 2392/81 = 11960/405 is the only value such as the Tzolk'in-lunar  commensuration period LCM(260,$T_i^0$) $<$ 18980 = 1 CR. The Copan Moon ratio is related to the Palenque formula by the relation LCM(4400,11960) = 1315600 = 44550 $\times$ 11960/405 = 44551 $\times$ 4400/149. The only other Moon ratio having similar properties is 383 lunations = 11310 days (11310/383 = LCM(260,11310)/766 = 22620/766 = 29.530026, $\varepsilon=0$) such as LCM(11310,11960) = 1040520 = 35235 $\times$ 11960/405 = 35236 $\times$ 11310/383. This value is far from the Moon synodic period compared to the Palenque formula and the Copan Moon ratio, maybe the reason why it is not present on Maya codices and monuments. 

Based on the hypothesis that Mayan priests-astronomers had knowledge of the planet synodic periods and the basic parameters for solar/lunar eclipse prediction, the model described above leads unequivocally to the Palenque and Copan lunar equations, providing a strong argument in favor of the use of this model by the Maya. The Mayan arithmetical model of astronomy can then be described as follows. The Calendar Round can be used to keep track of the solar year (Haab') and the synodic movement of Venus and Mars: LCM(260,365) = 18980 = 1CR, LCM(260,584) = 37960 = 2 CR and LCM(365,780) = 56940 = 3 CR, the length of the Dresden Codex Venus and Mars tables [\onlinecite{Bricker2011_163}]. The Tun-Haab'-Kawil wheel $\mathcal{Y}$ = LCM(360,365,3276) = 126 $\times$ 18980 = 2391480 days induces the movement of the wheels describing the synodic movement of Mercury, Jupiter, Saturn and the lunar months \{LCM($P_i$,$\mathcal{Y}$)/$\mathcal{Y}$, $i$ = $1\dots9$\} = \{29, 1, 1, 1, 19, 3, 59, 89, 37\} and gives rise to the calendar super-number $\mathcal{N}$. $\mathcal{N}$ can be expressed as:

\begin{eqnarray}\label{model}
\mathcal{N}  &=&   \mathcal{Y} \times 3 \times 19 \times 29 \times\ 37 \times 59 \times 89\\\nonumber
&\equiv& \alpha \times S_0
\end{eqnarray}
where $\alpha$ is an integer number of lunations: $\alpha$ = 26008014145502 corresponds to the Palenque formula $S_0$ = 2392/81 = 29.530864 days. The calendar super-number is a whole multiple of lunation given by the Palenque formula. The use of the Palenque formula is attested in several Classic period (200 to 900 AD) Maya sites and in Mayan codices up to the post-Classic period (1300 to 1521 AD). A calculation of the lunar series from recent excavations has shown that the Palenque formula was also used in Tikal (present-day Guatemala) on 9.16.15.0.0 or 17 February 766 AD, suggesting the widespread use of the Palenque formula throughout the Maya world [\onlinecite{Fuls2007}]. The Mayan arithmetical model of astronomy is equivalent to a clockwork mechanism in which the gears governing the planet synodic motions are driven by the synodic movement of the Moon, in the same sense as the contemporary representation of the Calendar Round (Fig. \ref{CR}). The Maya were aware of the imperfection of the model and attempted to improve it. Evidence of a X century AD Mayan astronomical innovation has been found in Chichen Itza (present-day Mexico) and in the Dresden Codex Venus table [\onlinecite{Aldana2016}]. 

\section{Astronomical derivation of the Mayan Calendar}

The Tzolk'in, the Tun, the Haab' and the 4-Kawil can be defined from the synodic periods of the Moon and the five planets observed with a naked-eye (Equ. \ref{astro}). For larger period of time, the Baktun (400 Tun = 144000 days) commensurates with the Haab' and the 4-Kawil forming the calendar grand cycle $\mathcal{GC}$ = LCM(365,3276,144000) = 400 $\times$ LCM(360,365,3276) = 956592000. The Euclidean division of $\mathcal{N}/37$ (LCM of 8 astronomical input parameters omitting for now the pentalunex) by $\mathcal{GC}$ gives:
\begin{eqnarray}\label{astro2}
\mathcal{N}/37&=&\mathcal{GC} \times \mathcal{Q}+ \mathcal{R}\\\nonumber
\mathcal{Q} &=& 21699 \\\nonumber
\mathcal{R} &=& 724618440\\\nonumber
&=& 101 \times126\times56940\\\nonumber
 &= & 126\times\sum_{i=0}^3\mathcal{X}_i.
\end{eqnarray}

If we note the Maya Aeon $\mathcal{A}$ = 13 $\times$ 73 $\times$ 144000 = 400 $\times$ $\mathcal{X}_0$ = 100 $\times$ $\mathcal{LR}$ = 136656000 such as $\mathcal{GC}$ = 7 $\times$ $\mathcal{A}$, the Euclidean division of $\mathcal{N}/37$ by $\mathcal{A}$ gives:
\begin{eqnarray}\label{astro1}
\mathcal{N}/37 &=&  \mathcal{A}  \times\mathcal{Q} +\mathcal{R}\\\nonumber
\mathcal{Q} &=& 151898\\\nonumber
\mathcal{R} &=& 41338440 \\\nonumber&=& 6\times121 \times 56940\\\nonumber &=& 121\times\mathcal{X}_0
\end{eqnarray}

The Maya Aeon such that $\mathcal{A}$ = LCM(260,365,144000) = 7200 $\times$ 18980 =  3600 $\times$ 37960 = 2400 $\times$ 56940 is commensurate to the 7200-day Katun, the Calendar Round, Venus and Mars synodic periods such as LCM(365,584) = 37960 and LCM(365,780) = 56940. We can rewrite Equ. \ref{astro2} and \ref{astro1} as:
\begin{eqnarray}\label{astro3}
\mathcal{N}/37 &-& 121 \times \mathcal{X}_0 = 151898\times \mathcal{A} \\\nonumber
\mathcal{N}/37 &-& 126 \times  \sum_{i=0}^3\mathcal{X}_i = 151893\times \mathcal{A} 
\end{eqnarray}

The subtraction of the two equations in Equ. \ref{astro3} can be expressed as a function of the Xultun numbers:
\begin{eqnarray}\label{astro4}
5\times\mathcal{A} &=&5  \times\mathcal{X}_0 + 95 \times 126 \times 56940\\\nonumber
5\times\mathcal{A} &=&5  \times\mathcal{X}_0 + \mathrm{LCM}(\sum_{i=1}^3\mathcal{X}_i ,\ \mathcal{X}_1 + 2 \mathcal{X}_2+ \mathcal{X}_3 )
\end{eqnarray}

The same results are obtained including the pentalunex with a factor 37. The grand cycle is such as $\mathcal{GC}$ = LCM(260,365,3276,$\mathcal{E}$) = $511\times\mathcal{E}$ where $\mathcal{E}$ = 1872000 days is the 13 Baktun Maya Era. The 5 Maya Aeon $5\mathcal{A}$ (Equ. \ref{astro4}) is such that 5 $\times$ $\mathcal{A}$ = 12000 $\times$ 56940 = 365 $\times$ $\mathcal{E}$. Equ. \ref{astro4} defines two important mythical dates $5\mathcal{X}_0$ and the Maya Era $\mathcal{E}$. We have: $5\times\mathcal{A} = 5\times\mathcal{X}_0 +570\times\mathcal{X}_1$ and $\mathcal{E}-5\times\mathcal{X}_0=1872000 - 1708200$ = 163800 = 10 $\times$ LCM(260,3276). $5\mathcal{X}_0$ has then the same properties as $5\mathcal{A}$ (same Tzolk'in, Haab', Kawil and direction-color) and the Maya Era $\mathcal{E}$ has the same Tzolk'in, Kawil and direction-color as $5\mathcal{A}$. The presence of the four Xultun numbers next to the Palenque formula on the walls of the small chamber excavated in Xultun provides an additional argument in favor of the Mayan arithmetical model of astronomy developed in section 2.

At this point, a question arises how the Maya, as early as the IX century AD, were able to compute tedious arithmetical calculations on such large numbers with up to 14 digits in decimal basis. Here is a possible method. From the prime factorizations of the canonic synodic periods $P_i$ (Table \ref{SP}), they may have listed all primes $p_i$ with their maximal order of multiplicity $\alpha_i$ and calculated the calendar super-number $\mathcal{N}$ (the LCM of the $P_i$'s) by multiplying each $p_i$'s $\alpha_i$ time. The Haab', the Tun and the 4-Kawil may have been obtained as described in section 2. From there, they may have proceeded to basic arithmetical calculations. The Euclidean division of $\mathcal{N}/37$ by $\mathcal{GC}$ = 7 $\times$ $\mathcal{A}$ = 400 $\times$ 7 $\times$ $\mathcal{X}_0$ (Equ. \ref{astro2}) is equivalent to a simplification of $\mathcal{N}/37$ by 7 $\times$ $\mathcal{X}_0$ = LCM(360,365,3276) = 2391480 and the Euclidean division of the product of the 5 left primes (3 $\times$ 19 $\times$ 29 $\times$ 59 $\times$ 89 = 8679903) by 400. The Euclidean division of $\mathcal{N}/37$ by $\mathcal{A}$ = $\mathcal{GC}/7$ = 400 $\times$ $\mathcal{X}_0$ (Equ. \ref{astro1}) is equivalent to a simplification of $\mathcal{N}/37$ by $\mathcal{X}_0$ = LCM(260,360,365) = 341640 and the Euclidean division of the product of the 6 left primes (3 $\times$ 7 $\times$ 19 $\times$ 29 $\times$ 59 $\times$ 89 = 60759321) by 400. It is to be noted that only 12 different prime factors $<100$ appear in the prime factorization of the calendar super-number which facilitates the operation (there are only 25 prime numbers lower than 100).

The Calendar Round and the Kawil-directions-colors wheel are initialized at the mythical origin of time as 4 Ahau 8 Cumku, 3 East-Red. This corresponds to a set of 4 indices that can be deduced from the calendar super-number as follows. We first create ordered lists of the Haab' and the Tzolk'in, assigning a unique set of 2 numbers for each day of the 18980-day Calendar Round [\onlinecite{Aveni2001_143,Aveni2001_147}]. For the Haab', the first day is 0 Pop (numbered 0) and the last day 4 Uayeb (numbered 364). For the Tzolk'in, the first day is 1 Imix (numbered 0) and the last day 13 Ahau (numbered 259). In this notation, the date of creation 4 Ahau 8 Cumku is equivalent to \{160;349\} and a date $D$ in the Calendar Round can be written as $D\equiv$ \{mod($D$ + 160,260);mod($D$ + 349,365)\}. The calendar super-number is such that: 
mod($\mathcal{N}$/13/37/73,260) = 160, mod($\mathcal{N}$/13/37/73,13) = 4,
mod($\mathcal{N}$/13/37/73,20) = 0 and mod($\mathcal{N}$/13/37/73,73) = 49. The Dresden Codex lunar equation 405 lunations = 11960 days is such that 405 $\times$ LCM(360,365,3276) = 405 $\times$ LCM(360,365,3276) = 405 $\times$ 2391480 = 51030 $\times$ 18980 = 80982 $\times$ 11960 + 4680. Starting from 4 Ahau 8 Cumku \{160;349\}, the date 80982 $\times$ 11960 = 968544720 days corresponds to 4 Ahau 8 Zip  \{160;49\}. The Calendar Round initialization at the origin of time may be related to the completion of a 11960-day lunar/solar eclipse cycle, 4680 days before the completion of the Tun-Haab'-Kawil commensuration cycle.The Kawil-directions-colors indices can be initialized  at the LCC origin of time as mod($\mathcal{N}$/37/32760,4) = 3 East-Red. That defines the position of the Calendar Round and the Kawil-directions-colors indices at the mythical date of creation, the Long Count date 13(0).0.0.0.0 4 Ahau 8 Cumku  \{160;349\}, 3 East-Red. 

\section{Mayan mythical dates and their numerological significance}
  
\begin{table}[h!]
\begin{center}
\begin{tabular}{l l l l l}
Date & $D$ & LCC date &  Cyclical date\\\hline
$\mathcal{I}_0$ & 0 & 13(0).0.0.0.0 &  \{160,349,3,0\}\\
$5\mathcal{X}_0$ & 1708200 & 11.17.5.0.0  & \{160,349,588,1\}\\
$\mathcal{E}$ & 1872000 & 13(0).0.0.0.0  & \{160,264,588,1\}\\
$5\mathcal{A}$ & \multicolumn{2}{l}{365$\times$13(0).0.0.0.0} &  \{160,349,588,1\}\\
$\mathcal{GC}$ & \multicolumn{2}{l}{511$\times$13(0).0.0.0.0} &  \{160,349,3,0\}\\
\end{tabular}

\caption{Important Mayan mythical dates: $\mathcal{I}_0$ (mythical date of creation), $5\mathcal{X}_0$ (date of the Itza prophecy), $\mathcal{E}$ (end of the 13 Baktun Era), $5\mathcal{A}$ (5 Maya Aeon) and $\mathcal{GC}$ = $7\mathcal{A}$ (Maya grand cycle). A date is defined as its linear time day $D$ and its cyclical equivalent given by the LCC date and a set of 4 integers \{$T$;$H$;$K$;$n$\} where $T$ is the Tzolk'in, $H$ the Haab', $K$ the Kawil and $n$ the direction-color calculated from the definitions given in the text.}\label{date}
\end{center}
\end{table}

\begin{figure}[h!]
\begin{center}
\includegraphics[width = 70mm]{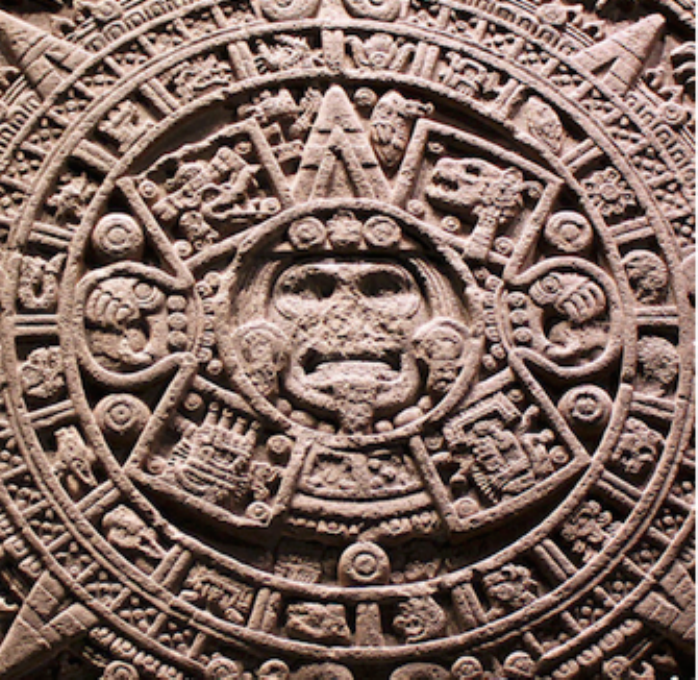}
\caption{Aztec Calendar Stone representing the Five Creations of the Mayan/Aztec civilization. The stone, probably carved during the post-Classic period (1300 to 1521 AD),  was rediscovered in 1790 during the construction of the cathedral of Mexico City, Mexico.}\label{SunStone}
\end{center}
\end{figure}

The Mayan cyclical concept of time implied that current and future events were pre-determined from past events in an ever-repeating cosmological grand cycle. For divinatory purposes, Mayan priests-astronomers developed an elaborate numerology based on the Mayan arithmetical model of astronomy. This numerology allowed them to discriminate between two dates over very long period of time: the combination of the Tzolk'in, the Haab' and the Kawil-directions-colors differentiates two dates every $\mathcal{X}_1$ = LCM(260,365,3276) =1195740 days $\approx$ 3273 years and the 13 Baktun Maya Era every 1872000 days $\approx$ 5125 years. Table \ref{date} gives the full Mayan Calendar dates of important mythical dates deduced from Equ. \ref{astro4}. The previous Maya Era is characterized by two important dates: the mythical date $5\mathcal{X}_0$ = 11.17.5.0.0 4 Ahau 8 Cumku, 588 South-Yellow and the end of the 13 Baktun Era $\mathcal{E}$ = 13(0).0.0.0.0 4 Ahau 3 Kankin, 588 South-Yellow, which are defined by their equivalent properties compared to the 5 Maya Aeon cycle $5\mathcal{A}$ (4 Ahau 8 Cumku, 588 South-Yellow). The 5 Maya Aeon cycle is consistent with the Aztec myth of the Five Suns as represented by the Calendar Stone (Fig. \ref{SunStone}): the central solar deity corresponding to the Fifth Sun or current Creation is surrounded by four squares representing the four previous Creations. Indeed, the Maya (100 BC to 1521 AD) and the Aztec (1300 to 1521 AD), both originating from the ancient Olmec civilization, shared many similarities as regards their myths and beliefs. The Aztec calendar, adapted from the Mayan Calendar, was a combination of the 365-day Xiuhpohualli year cycle and the 260-day Tonalpohualli ritual cycle coinciding every 52 years = 1 Calendar Round. The date $5\mathcal{X}_0$ or 3 July 1564 AD (GMT correlation) may be related to the Itza prophecy predicting an intense cultural change of the Aztec civilization during the XVI century AD [\onlinecite{VanStone,Stuart2012_19}]. The occurrence of a series of omens prior to the arrival of the Spaniards was interpreted as the upcoming fulfillment of the prophecy, thought to be the annunciation of the Spanish conquest of Mexico (from February 1519 to 13 August 1521).

\begin{figure}[h!]
\begin{center}
\includegraphics[width = 90mm]{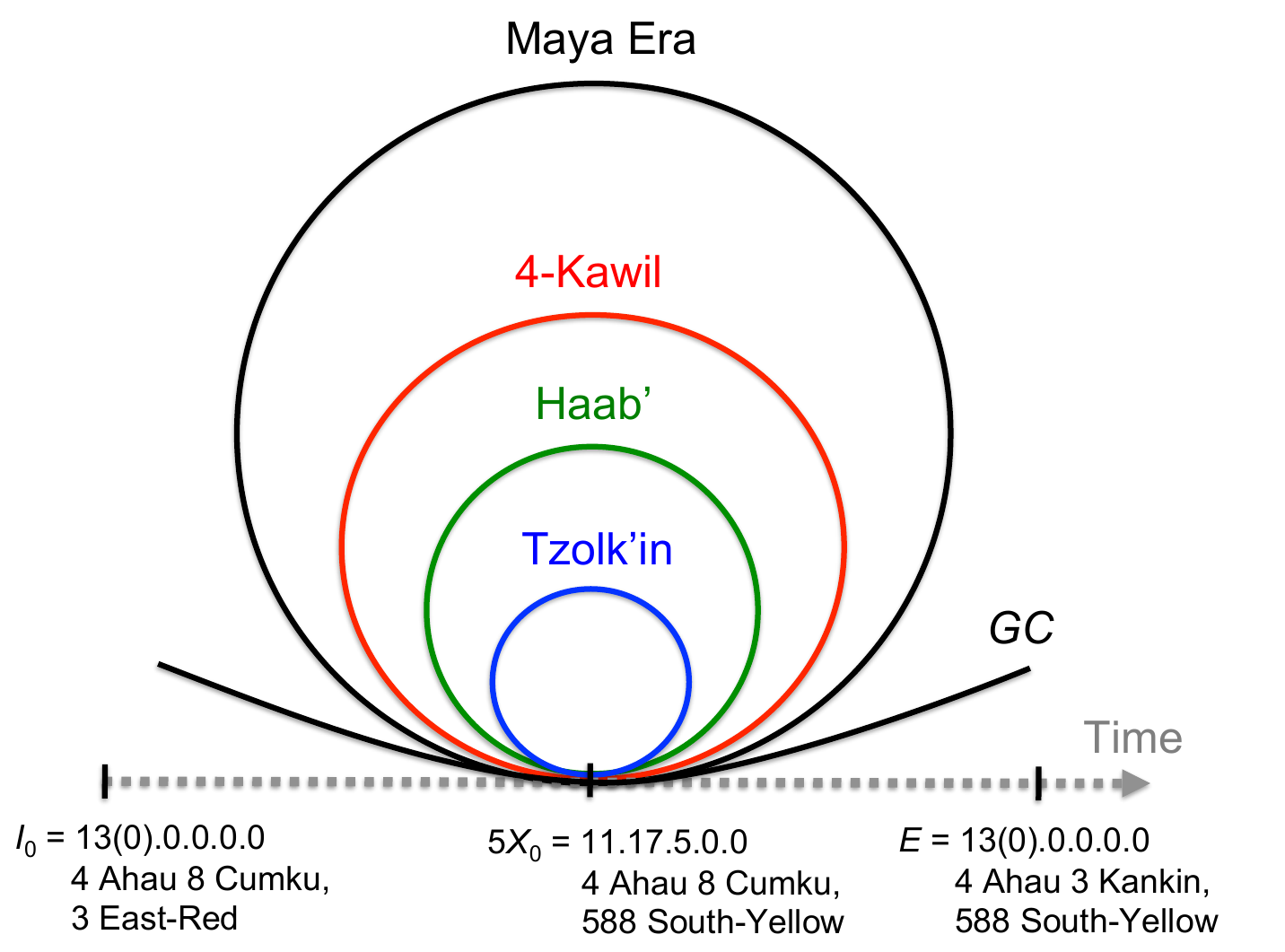}\\[.5cm]
\caption{Mayan cyclical $vs$ linear concept of time, with the  260-day Tzolk'in, the 365-day Haab', the 3276-day 4-Kawil, the Maya Era of 13 Baktun (5125 years) and the grand cycle $\mathcal{GC}$ of 511 Maya Eras. The mythical date of creation $\mathcal{I}_0$ (11 August 3114 BC), the date of the Itza prophecy $5\mathcal{X}_0$ (3 July 1564 AD) and the end of the previous Maya Era $\mathcal{E}$ (21 December 2012) are also represented.\\[-1cm]}\label{ZZZ}
\end{center}
\end{figure}

Fig. \ref{ZZZ} represents the Mayan cyclical concept of time, with a grand cycle $\mathcal{GC}$ defined as the commensuration of the Tzolk'in, the Haab', the Kawil-directions-colors and the Maya Era. The three important dates of the previous Maya Era are represented: the mythical date of creation 13(0).0.0.0.0 4 Ahau 8 Cumku, 3 East-Red (11 August 3114 BC), the date corresponding to the Itza prophecy 11.17.5.0.0 4 Ahau 8 Cumku, 588 South-Yellow (3 July 1564 AD) and the end of the Maya Era 13(0).0.0.0.0 4 Ahau 3 Kankin, 588 South-Yellow (21 December 2012). Mayan priests-astronomers may have chosen the end of the 13 Baktun to coincide with the winter solstice of 2012 as suggested earlier [\onlinecite{Edmonson}]. The presence of tropical year-related Long Count numbers in Maya inscriptions is discussed in Ref. [\onlinecite{Grofe}]. If we assume their knowledge of the calendar year drift cycle 1508 Haab' = 1507 tropical year corresponding to a very good approximation of the tropical year of 365.2422 days, starting from the GMT correlation date of origin 11 August 3114 BC, we arrive at 21 December 2012, 13 Baktun = 1872000 days later.\\[1cm]

\section{Conclusion}
This study presents a complete description of the Mayan cyclical concept of time, characterized by a set of calendar cycles derived from an arithmetical model of naked-eye astronomy. Based on a realistic hypothesis that Mayan priets-astronomers had a good knowledge of naked-eye astronomy, this model is derived from a set of 9 integer input parameters, the solar year (Haab'), the three lunar months (the pentalunex and the two lunar semesters) and the synodic periods of Mercury, Venus, Mars, Jupiter and Saturn and was used to determine the Mayan lunar equation from astronomical observations. The Palenque formula and the Copan Moon ratio are solutions of the model and were certainly derived in or prior to the Maya Classic period (200 to 900 AD). The results are expressed as a function of the IX century AD Xultun numbers, four enigmatic Long Count numbers deciphered on the walls of a small room in the Maya ruins of Xultun (present-day Guatemala), providing additional argument in favor of the model. The calendar super-number derived from the model leads to the Mayan Calendar cycles: the 3276-day Kawil-directions-colors, combination of the 4 directions-colors and the 819-day Kawil, the 18980-day Calendar Round, combination of the 260-day Tzolk'in and the 365-day Haab', and the 1872000-day Long Count Calendar formed by the 360-day Tun, the 7200-day Katun and the 144000-day Baktun. The Mayan cyclical concept of time is explained, in particular the position of the Calendar Round at the mythical date of creation 13(0).0.0.0.0 4 Ahau 8 Cumku and the numerological significance of the 13 Baktun Maya Era. The Mayan Calendar, derived from the Mayan arithmetical model of astronomy, shows the high proficiency of Mayan mathematics applied to numerological timekeeping for divinatory purposes.

\onecolumngrid

\vline

$^*$ e-mail: \href{mailto:thomas.chanier@gmail.com}{thomas.chanier@gmail.com}

\twocolumngrid




\end{document}